\documentclass{article}

\usepackage{graphicx}
\usepackage{amsmath}
\usepackage[dvipsnames]{xcolor}
\usepackage[margin=1.5in]{geometry}

\newcommand{\fn}[0]{\nabla_{\bX}}
\newcommand{\sn}[0]{\nabla_{\bx}}
\newcommand{\dd}{\partial}
\newcommand{\del}{\nabla}
\renewcommand{\d}{{\rm d}}
\newcommand{\beq}{\begin{equation}}
\newcommand{\eeq}{\end{equation}}
\newcommand{\bx}{{\bf x}}
\newcommand{\bX}{{\bf X}}
\newcommand{\bn}{{\bf n}}
\newcommand{\bQ}{{\bf Q}}
\newcommand{\bV}{{\bf V}}
\newcommand{\ra}{\rightarrow}
\newcommand{\ii}{{\mathrm i}}
\newcommand{\ee}{{\mathrm e}}

\title{Integral Constraints in Multiple Scales Problems with a Slowly Varying Microstructure}

\author{A. Kent$^1$, S. L. Waters$^1$, J. Oliver$^1$, S. J. Chapman$^1$}

\date{$^1$ Mathematical Institute, University of Oxford, Woodstock Road, Oxford OX2 6GG.}

\begin{document}

\label{firstpage}
\maketitle

\begin{abstract}
  Asymptotic homogenisation is considered for problems with integral constraints imposed on a slowly-varying microstructure; an insulator with an array of perfectly dielectric inclusions of slowly varying size serves as a paradigm. Although it is well-known how to handle each of these effects (integral constraints, slowly-varying microstructure) independently within multiple scales analysis,
additional care is needed when they are combined. Using the flux transport theorem, the multiple scales form of an integral constraint on a slowly varying domain is identified. The proposed form is applied to obtain a homogenised model for the electric potential in a dielectric composite, where the microstructure slowly varies and the integral constraint arises due to a statement of charge conservation. A comparison with multiple scales analysis of the problem with established approaches provides validation that the proposed form results in the correct homogenised model. 

\vspace{0.5em}

\vspace{0.5em}
\end{abstract}

\section{Introduction}
Homogenisation via multiscale asymptotics is one coarse graining method that can be used to derive the effective properties of composite media \cite{Hornung1996}. Typically used for periodic microstructure, example applications of the technique include modelling flow in porous media, wave propagation in poroelastic materials, filtration and decontamination processes \cite{Burridge1981, Dalwadi2015, Hornung1996,Luckins2020}.

The result of the homogenisation process is the reduction of a problem posed on a complicated domain, or with rapidly varying coefficients, to two simpler problems: one `cell problem' describing the microscale variation; and a  second  `homogenised model' describing the macroscale variation of variables across the whole domain. The technique can be extended to problems with a slowly varying geometry, albeit at the cost of having a  cell problem which varies with the macroscale \cite{Dalwadi2017}. A mapping depending on the slow spatial scale can be applied to transform a heterogeneous microstructure to an exactly periodic reference configuration \cite{Hornung1992, Peter2007}. Standard homogenisation can be performed in this reference configuration before inverting the mapping to obtain the homogenised equations featuring spatially-dependent coefficients which reflect microstructural variation. A similar approach can be used to treat microstructures with temporal and spatiotemporal variations. Examples of problems using a prescribed mapping include \cite{Bruna2015, Collis2017, Richardson2011} and the mapping can be coupled to the macroscale variables \cite{Peter2009}. When the domain is locally periodic and the unit cell has fixed size, transformation to a reference configuration is no longer required as the slow variable features as a parameter in the microscale problem \cite{Dalwadi2017}. Examples of this approach are found in \cite{Dalwadi2015, Fatima2011, Ray2015}. When considering problems with a slowly varying domain, care must be taken in converting Neumann and Robin boundary conditions on microscopic inclusions into multiple scales form. Typically, a level set function is introduced to define the boundary of the inclusion, with the expansion of the normal to the boundary derived by writing the gradient of the level set function in multiple scales form \cite{Bruna2015, Fatima2011, Richardson2011}.  

A second extension of the standard method allows for homogenisation of problems featuring integral constraints \cite{Chapman2015}. These constraints generally appear as conservation conditions, for example,  of charge or momentum, with applications in modelling nematic crystals, radiation in porous media and bubbly liquids \cite{Bennett2018, Chapman2015, Rooney2021}. Unlike standard multiple scales, where the macroscale coordinate can be assumed to take a constant value within a given unit cell, it is crucial to account for the small variation in macroscale coordinate along the integration path, since this variation causes a change in flux which  affects the parameters in the homogenised model.  

In the present work we aim to combine these two extensions, developing an understanding of how to write integral constraints on a slowly varying domain in multiple scales form. Although this seems like a routine task, we will see that in fact that the answer is not obvious a priori.
We use as a paradigm the problem of the electric potential in an insulator interspersed with a periodic array of perfectly dielectric inclusions of slowly varying size. This problem has the advantage that the perfectly dielectric limit can also be taken after a standard homogenisation procedure, so that we know what the homogenised model should be. Not all integral constraint problems can be recast in this way.

\section{Paradigm Problem}\label{sec:2d}
\begin{figure}
    \centering
    \includegraphics[width = 0.6\textwidth]{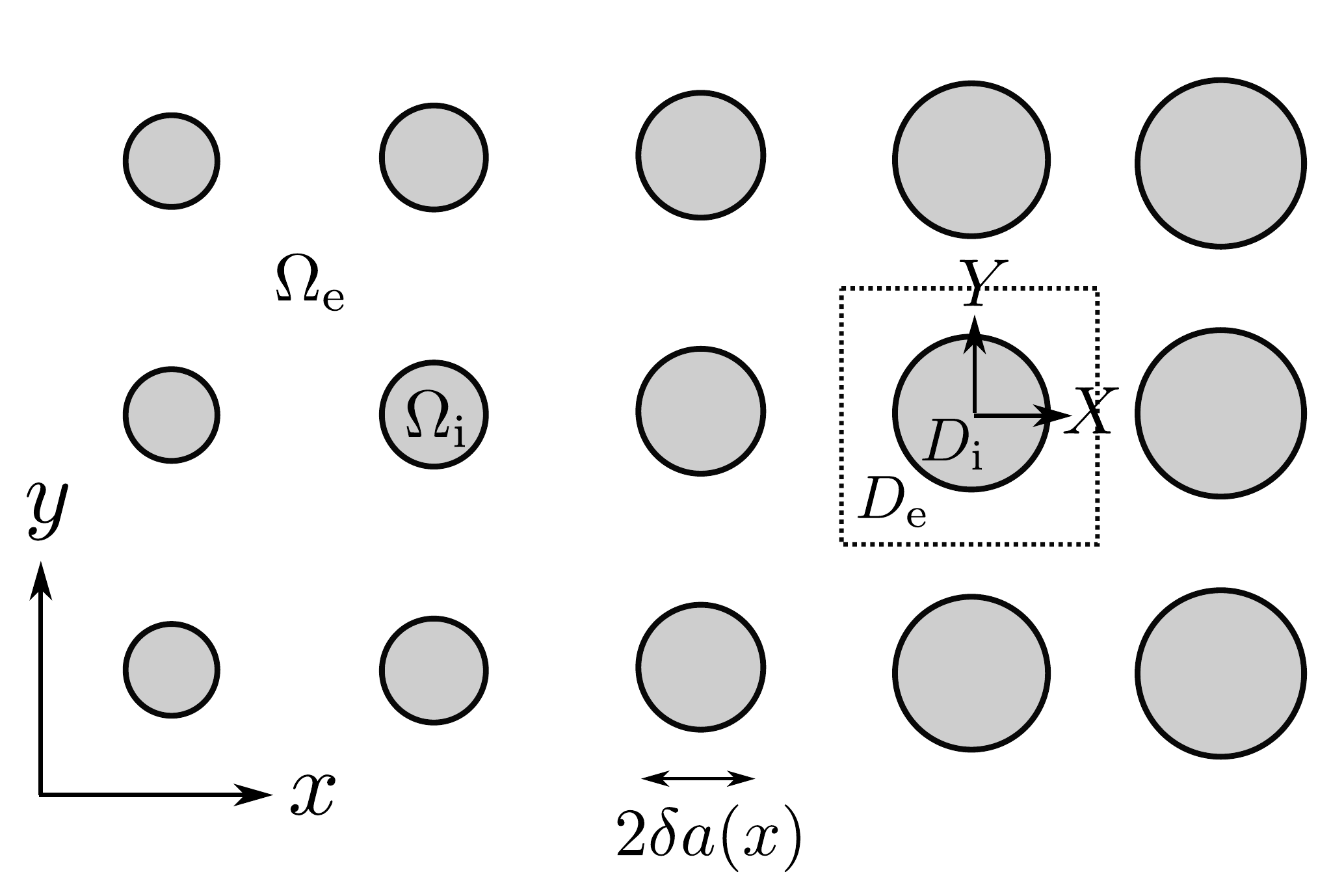}
    \caption{Schematic of the 2D composite. \textit{Perfectly dielectric inclusions shown in grey lie on a periodic array within an insulator. The inclusions have a radius $a(x)$ which slowly varies across the domain. } }
    \label{fig:pot2d}
\end{figure}

We consider the electric potential $\phi$ in a
dielectric material, which satisfies Poisson's equation
\begin{equation}
    \nabla \cdot ( \varepsilon \nabla \phi) = -\rho \label{eq:2dgauss},
\end{equation}
where $\varepsilon$ is the permittivity and $\rho$ is the charge density (which we suppose is given). We consider a material composite comprising  an insulator $\Omega_\ee$ of constant permittivity $\varepsilon_\ee$ with an array of inclusions $\Omega_\ii$ of constant permittivity $\varepsilon_\ii$.
At the boundary between the two regions
\begin{align}
    \left[ \bn \cdot( \varepsilon \nabla \phi ) \right]^\ee_\ii &= 0, \label{eq:2decont} \\
    \left[ \phi \right]^\ee_\ii &= 0, \label{eq:2dcont}
\end{align}
where $\bn$ is the (outward-facing) normal to the boundary of the inclusion, and $[\cdot]^\ee_\ii$ represents the jump in the enclosed quantity across the interface.
We suppose that the centres of the inclusions are arranged on a regular cubic lattice of side $\delta$, and that the radius of each inclusion $\delta a(\bx)$ varies slowly with (macroscopic) position (see Fig.~\ref{fig:pot2d}).

In the limit $\varepsilon_\ii \to \infty$ the inclusions are perfectly dielectric and the model becomes 
\begin{align}
  \nabla \cdot ( \varepsilon_\ee \nabla \phi) & = -\rho \qquad \mbox{ in }\Omega_\ee,\label{lim:eqe}\\
  \nabla \phi & = 0\qquad \mbox{ in }\Omega_\ii,\label{lim:eqi}
\end{align}
with boundary condition
\begin{align}
    \left[ \phi \right]^\ee_\ii &= 0. \label{lim:2dcont}
\end{align}
The potential $\phi$ is constant on each inclusion, but may take different values on different inclusions. To close the problem we need to integrate \eqref{eq:2dgauss} over each inclusion and use \eqref{eq:2decont} to give the integral constraint
\beq
\int_{\dd \Omega_\ii} \left.\varepsilon_\ee\, \bn\cdot \del \phi \right|_\ee\, \d S = - \int_{\Omega_\ii} \rho \, \d \bx, \label{lim:int}
\eeq
where $\Omega_\ii$ is an individual inclusion. 

We will approach the limit $\varepsilon_\ii \to \infty$ in two different ways.
We will first homogenise (\ref{eq:2dgauss})-(\ref{eq:2dcont}) following \cite{Bruna2015}, before  taking the limit $\varepsilon_\ii \to \infty$ in the homogenised model.
We will then homogenise \eqref{lim:eqe}-\eqref{lim:int} directly, which will require us to determine how to cast \eqref{lim:int} in  multiple-scales form when the domain $\Omega_\ii$ is a function of (slow) position.

\subsection{Standard Multiple Scales}\label{sec:2dstandard}
We introduce the fast scale $\bX = \bx/\delta$, and
suppose that $\phi = \phi(\bx,\bX)$, treating $\bx$ and $\bX$ as independent,
with derivatives transforming according to the chain rule
\beq
\del \ra \sn + \frac{1}{\delta}\fn.\label{chain}
\eeq
We remove the indeterminancy that this introduces by imposing that $\phi$ is 
${\bf 1}$-periodic in  $\bX$.
To describe the inclusions we introduce the function 
\begin{equation}
    h(\bx,\bX) = |\bX - \lfloor \bX \rfloor| - a(\bx),\label{eq:levelset}
\end{equation}
where $ \lfloor \bX \rfloor$ represents the integer part of each component of $\bX$. This function is ${\bf 1}$-periodic in $\bX$, and the level set $h = 0$ defines the boundary of the inclusion. The normal to the inclusion can then be written in multiple scales form as
\begin{equation}
    \bn =\frac{\del h}{|\del h|} =  \frac{\fn h + \delta \sn h }{|\fn h + \delta \sn h |}  = \bn_0 + \delta \bn_1 + O(\delta^2 ), \label{eq:ndirect}
\end{equation}
where
\beq
\bn_0 = \frac{\fn h}{|\fn h|}, \qquad
\bn_1 = \frac{\sn h}{|\fn h|} - \frac{(\sn h \cdot \fn h) \fn h}{|\fn h|^3}.
\label{eq:normal}
\eeq
Within a given unit cell $D$ we denote the region occupied by the inclusion as $D_\ii(\bx)$ and that occupied by the insulator as $D_\ee(\bx)$, that is
\[
  D_\ii =  \{ \bX \in D: h(\bx,\bX)<0\}, \qquad
  D_\ee  =  \{ \bX \in D: h(\bx,\bX)>0\}.
\]
 Substituting (\ref{chain}) and \eqref{eq:ndirect} into
 (\ref{eq:2dgauss})-(\ref{eq:2dcont}), expanding
\begin{equation}
    \phi  \sim \phi_0(\bx, \bX)  +  \delta \phi_1(\bx, \bX) + \cdots,\label{eq:2dphimsexp}
\end{equation}
and equating coefficients of $\delta$, we find that at leading-order 
\begin{align}
\fn  \cdot ( \varepsilon \fn \phi_0) &= 0,\\
 \left[ \bn_0 \cdot( \varepsilon \fn  \phi_0 )  \right]^\ee_\ii &= 0, \\
 \left[ \phi_0 \right]^\ee_\ii &= 0,
\end{align}
with $\phi_0$ ${\bf 1}$-periodic in $\bX$. Thus, $\phi_0$ is constant on the fast scale, so that $\phi_0 = \phi_0(\bx)$. At next order, we find
\begin{align}
    \fn \cdot \left( \varepsilon \fn \phi_1 \right) &= 0, \\
     \left[ \bn_0 \cdot \left( \varepsilon(  \fn \phi_1 +  \sn  \phi_0)  \right)  \right]^\ee_\ii &= 0,\label{2.17} \\
 \left[ \phi_1 \right]^\ee_\ii &= 0,
\end{align}
with $\phi_1$ ${\bf 1}$-periodic in $\bX$, where we have used the fast scale independence of the leading-order potential. The solution is 
\begin{equation}
    \phi_1 = {\bf \Psi}\cdot \sn \phi_0 + \overline{\phi}_1,\label{eq:2dphi1}
\end{equation}
where 
$\overline{\phi}_1$ is independent of $\bX$ and ${\bf \Psi}$ satisfies the cell problem
\begin{align}
    \fn \cdot \left( \varepsilon \fn {\bf \Psi} \right) &= 0,  \label{eq:potcell2d}\\
    \left[ \bn_0  \cdot \left( \varepsilon (\fn {\bf \Psi} + I   )\right) \right]^\ee_\ii &= 0,\label{eq:bcncell} \\
    \left[ {\bf \Psi} \right]^\ee_\ii &= 0,\label{eq:potcellcont}
\end{align}
with ${\bf \Psi}$ ${\bf 1}$-periodic in $\bX$, where $I$ is the identity matrix, and uniqueness is achieved by imposing zero mean, for example. Finally, equating coefficients of $\delta^2$, we find
\begin{align}
  \fn \cdot \left( \varepsilon ( \fn \phi_2  + \sn \phi_1   )\right) + \sn \cdot \left( \varepsilon ( \fn \phi_1  + \sn \phi_0   )\right) &= -\rho,  \label{eq:2dod2gauss} \\
     \left[ \bn_0 \cdot \left( \varepsilon(  \fn \phi_2 +  \sn  \phi_1)  \right) \right]^\ee_\ii  +   \left[ \bn_1 \cdot \left( \varepsilon(  \fn \phi_1 +  \sn  \phi_0)  \right)  \right]^\ee_\ii &= 0, \label{eq:2dod2jump}\\
 \left[ \phi_2 \right]^\ee_\ii &= 0.\label{eq:2dod2cont}
\end{align}
Integrating (\ref{eq:2dod2gauss}) over the unit cell $D$, applying the divergence theorem to terms involving the fast divergence, and using (\ref{eq:2dod2jump}) we find
\begin{multline}
    \int_{\partial D_\ee} \varepsilon_\ee \left( \fn \phi_1  +  \sn \phi_0 \right) \cdot \bn_1  \d S_\bX -  \int_{\partial D_\ii}  \varepsilon_\ii \left( \fn \phi_1  +  \sn \phi_0 \right) \cdot \bn_1  \d S_\bX  \\ +\int_{D_\ee} \sn \cdot \left( \varepsilon_\ee ( \fn \phi_1  + \sn \phi_0   )\right)\d \bX +  \int_{D_\ii } \sn \cdot \left( \varepsilon_\ii ( \fn \phi_1  + \sn \phi_0   )\right) \d \bX = -\rho_{{\rm eff}},
\end{multline}
where $\partial D_\ii$ and $\partial D_\ee$ denote the interior and exterior of the inclusion boundary in the unit cell respectively, and the effective charge  is given by
\begin{equation}
    \rho_{{\rm eff}} = \int_{D} \rho \, \d \bX. \label{eq:rhoeff}
\end{equation}
Taking the slow divergence outside the integral using the Reynolds transport theorem, we find
\begin{multline}
    \int_{\partial D_\ee} \varepsilon_\ee \left( \fn \phi_1  +  \sn \phi_0 \right) \cdot \bn_1  \,\d S_\bX -  \int_{\partial D_\ii}  \varepsilon_\ii \left( \fn \phi_1  +  \sn \phi_0 \right) \cdot \bn_1  \,\d S_\bX  \\
       + \int_{\partial D_\ee} \varepsilon_\ee \left( \fn \phi_1  +  \sn \phi_0 \right) \cdot \bV   \cdot \bn_0 \, \d S_\bX -  \int_{\partial D_\ii}  \varepsilon_\ii \left( \fn \phi_1  +  \sn \phi_0 \right)\cdot \bV \cdot \bn_0  \, \d S_\bX \\
 +\sn \cdot   \int_{D_\ee}  \left( \varepsilon_\ee ( \fn \phi_1  + \sn \phi_0   )\right)\, \d \bX +   \sn \cdot \int_{D_\ii }  \left( \varepsilon_\ii ( \fn \phi_1  + \sn \phi_0   )\right) \, \d \bX = -\rho_{{\rm eff}},\label{2.28}
\end{multline}
where the matrix $\bV$ is the ``velocity'' of the boundary, i.e. the derivative of   position on the boundary with respect to $\bx$.
Differentiating the equation $h=0$ with respect to $\bx$ gives
\[\bV\cdot  \fn h + \sn h =0,\]
so that
\[ \bV \cdot \bn_0 = -\frac{ \sn h}{|\fn h|},
\qquad \bV \cdot \bn_0 + \bn_1 =  - \frac{(\sn h \cdot \fn h) \fn h}{|\fn h|^3} = - \frac{(\sn h \cdot \fn h)}{|\fn h|^2}\bn_0.\]
Thus, using (\ref{2.17}), the surface integrals cancel in (\ref{2.28}), leaving
\begin{equation}
     \sn \cdot   \int_{D_\ee}  \left( \varepsilon_\ee ( \fn \phi_1  + \sn \phi_0   )\right)\, \d \bX +   \sn \cdot \int_{D_\ii }  \left( \varepsilon_\ii ( \fn \phi_1  + \sn \phi_0   )\right)\, \d \bX = -\rho_{{\rm eff}}.
\end{equation}
Substituting (\ref{eq:2dphi1}), gives, finally, the homogenised problem 
\begin{equation}
  \sn \cdot \left( {\boldsymbol \varepsilon}_{{\rm eff}}  \sn \phi_0 \right)
  = -\rho_{{\rm eff}}, \label{eq:homprob}
\end{equation}
where the effective permittivity ${\boldsymbol \varepsilon}_{{\rm eff}} $ is given by
\begin{equation}
    {\boldsymbol \varepsilon}_{{\rm eff}}  = \int_{D} \varepsilon \left(  I + \fn {\bf \Psi} \right) \, \d \bX.\label{2.31}
\end{equation}

\subsubsection{The limit $\varepsilon_\ii \to \infty$}
As  $\varepsilon_\ii \to \infty$ in the cell problem   (\ref{eq:potcell2d})-(\ref{eq:potcellcont}) we find
\begin{align}
    \fn^2 {\bf \Psi}    &= 0   \quad \text{in} \quad D, \label{2.32}\\
   \bn_0 \cdot \left( \fn {\bf \Psi} + I \right) &= 0  \quad \text{on} \quad \partial D_\ii, \label{eq:inflimbc} \\
    \left[ {\bf \Psi}\right]^\ee_\ii &= 0.
\end{align}
Thus ${\bf \Psi} = -\bX + $ constant in $D_\ii$, where the constant must be chosen so that ${\bf \Psi}$ has zero mean. In the effective permittivity (\ref{2.31}) this gives zero times infinity in the inclusion, so we must manipulate this expression into something more suitable before we take the limit.
Switching to index notation, using \eqref{eq:bcncell}, the divergence theorem,   and (\ref{2.32}), we find 
\begin{eqnarray}
     {{\varepsilon}_{{\rm eff}}}_{ij}  & =& \int_{D} \varepsilon \delta_{ij} \, \d \bX + \int_{D} \varepsilon \frac{\dd}{\partial X_k} \left( X_j \frac{\dd \Psi_i}{\dd X_k} \right)\, \d \bX\nonumber\\
     & = &  \int_{D_\ii} \varepsilon_\ii  \delta_{ij} \, \d \bX +  \int_{D_\ee} \varepsilon_\ee \delta_{ij} \,\d \bX + \int_{\partial D_\ii} \varepsilon_\ii   X_j \frac{\dd\Psi_i }{\dd X_k } n_k \,\d S_\bX\nonumber \\&& \mbox{ }- \int_{ \partial D_\ee} \varepsilon_\ee   X_j \frac{\dd\Psi_i }{\dd X_k } n_k \,\d S_\bX + \int_{ \partial D} \varepsilon_\ee   X_j\frac{\dd\Psi_i }{\dd X_k } n_k \,\d S_\bX \nonumber\\
     & = & \int_{D_\ii} \varepsilon_\ii  \delta_{ij} \, \d \bX +  \int_{D_\ee} \varepsilon_\ee \delta_{ij}\, \d \bX - \int_{\partial D_\ii} \varepsilon_\ii  X_j  n_i  \,\d S_\bX \nonumber \\&& \mbox{ } + \int_{ \partial D_\ee} \varepsilon_\ee   X_j n_i  \,\d S_\bX + \int_{ \partial D} \varepsilon_\ee   X_j \frac{\dd\Psi_i }{\dd X_k } n_k  \,\d S_\bX\nonumber \\
     &=& \varepsilon_\ee \left( \delta_{ij} + \int_{\partial D} X_j \frac{\dd\Psi_i }{\dd X_k } n_k\,\d S_\bX  \right). \label{eq:permlim}
\end{eqnarray}
We can now safely take the limit $\varepsilon_\ii \to \infty$.

\subsection{Multiple Scales with Integral Constraints}\label{sec:2dintegrals}
We now apply the method of multiple scales directly to the problem
\eqref{lim:eqe}-\eqref{lim:int}, hoping to retrieve (\ref{eq:homprob}) with (\ref{eq:permlim}).
 Substituting (\ref{chain}) into
\eqref{lim:eqe}-\eqref{lim:2dcont}, expanding as in (\ref{eq:2dphimsexp}), and equating coefficients of $\delta$ we find that at  leading-order 
\begin{align}
    \fn\cdot ( \varepsilon_\ee \fn \phi_0) & = 0 \quad \text{in} \quad D_\ee, \\
    \fn \phi_0 &=\textbf{0} \quad \text{in} \quad D_\ii,\\
    \left[ \phi_0  \right]^\ee_\ii &= 0 ,
\end{align} 
with $\phi_0$ ${\bf 1}$-periodic in $\bX$. Thus, as before, $\phi_0 = \phi_0(\bx)$. At first-order we find
\begin{align}
    \fn\cdot ( \varepsilon_\ee \fn \phi_1 ) & = 0 \quad \text{in} \quad D_\ee, \label{eq:2dodgauss} \\
    \fn \phi_1 + \sn \phi_0 &= {\bf 0} \quad \text{in} \quad D_\ii,\label{eq:odinc}\\
    \left[ \phi_1  \right]^\ee_\ii &= 0 ,
\end{align}
with $\phi_1$ ${\bf 1}$-periodic in $\bX$. As in Section \ref{sec:2dstandard}, the solution is $\phi_1 = {\bf \Psi}\cdot \sn \phi_0 + \overline{\phi}_1$ where
$\overline{\phi}_1$ is independent of $\bX$ and
\begin{align}
    \fn \cdot ( \varepsilon_\ee \fn {\bf \Psi} ) &= {\bf 0} \quad \text{in} \quad D_\ee,\label{eq:cellint} \\
    \fn {\bf \Psi} + I &= 0 \quad \text{in} \quad D_\ii, \\
    \left[ {\bf \Psi} \right]^\ee_\ii &= {\bf 0} ,\label{eq:cellintcont}
\end{align}
with  ${\bf \Psi}$ ${\bf 1}$-periodic in $\bX$, and we impose 
\begin{equation}
    \int_{D} {\bf \Psi} \, \d \bX = {\bf 0}.
\end{equation}
Equating coefficients of $\delta^2$ we find
\begin{align}
    \fn \cdot \left( \varepsilon_\ee (\fn \phi_2 + \sn \phi_1 )\right) +   \sn \cdot \left( \varepsilon_\ee (\fn \phi_1 + \sn \phi_0 )\right) &= - \rho  \quad \text{in} \quad D_\ee, \label{eq:od2gaussint}\\
    \fn \phi_2 + \sn \phi_1 &= {\bf 0} \quad \text{in} \quad D_\ii,  \label{eq:od2gaussinc}\\
    \left[ \phi_2 \right]^\ee_\ii &= 0. 
\end{align}
Integrating (\ref{eq:od2gaussint}) over the exterior region and applying the divergence theorem to the first term, gives
\begin{equation}
  -  \int_{\partial D_\ee}   \varepsilon_\ee (\fn \phi_2 + \sn \phi_1 )\cdot \bn_0 \, \d S_\bX + \int_{D_\ee} \sn \cdot \left( \varepsilon_\ee (\fn \phi_1 + \sn \phi_0 )\right)\, \d \bX = -\int_{D_\ee} \rho \, \d \bX,\label{2.49}
\end{equation}
where the integral over the exterior boundary of the unit cell vanishes due to periodicity.
To evaluate the surface integral in \eqref{2.49} we need to use the integral constraint \eqref{lim:int}.

\subsubsection{Dealing with the integral}
As discussed in \cite{Chapman2015}, it seems natural to write \eqref{lim:int} in multiple scales form as
\[
\delta^2\int_{D_\ee} \varepsilon_\ee\, \bn\cdot\left( \sn \phi + \frac{1}{\delta}\fn \phi\right) \, \d S_\bX = - \delta^3\int_{\Omega_\ii} \rho \, \d \bX,
 \]
 but this is incorrect, as it neglects the small variation in the slow variation $\bx$ around the boundary of the inclusion, which turns out to be crucial.
 Writing $\bQ = \del \phi$, the approach taken in \cite{Chapman2015} was to recognise that on the interface $\bx = \hat{\bx} + \delta \bX$, where $\hat{\bx}$ is the position of the bottom left corner of the unit cell and $\bX \in [0,1]^3$, expanding
 \beq
 \bQ(\bx,\bX) =\bQ(\hat{\bx}+\delta \bX,\bX)  = \bQ(\hat{\bx},\bX)+\delta \bX \cdot \sn \bQ(\hat{\bx},\bX)+\cdots \label{eq:fastslow}
 \eeq
 in the integrand of \eqref{lim:int}. But how should we proceed when the domain and the normal, as well as the integrand, depend on the slow variable $\bx$?

\vspace{1em}

\noindent \textit{i.) The naive approach}
\vspace{0.5em}

\noindent An initial attempt to write the integral constraint on a slowly varying domain in multiple scales form may be to combine the normal expansion \eqref{eq:ndirect} with the expansion of the integrand given in \eqref{eq:fastslow}, writing
\begin{equation}
    \int_{\partial \Omega_\ii} \bQ \cdot \bn dS \to \delta^2 \int_{\partial \Omega_\ii} \left(\bQ_0 + \delta( \bQ_1 + \bX  \cdot \sn \bQ ) + ...    \right) \cdot (\bn_0 + \delta \bn_1 + ...) \, \d S_\bX. \label{eq:naive}
\end{equation}
We will now highlight the issues that arise if the form \eqref{eq:naive} is used. Writing \eqref{eq:od2gaussint}-\eqref{eq:od2gaussinc} in terms of the flux $\bQ$, and integrating over the insulating region, we find
\beq
\int_{D_\ee} \fn \cdot \bQ_1 \, \d \bX + \int_{D_e} \sn \cdot \bQ_0 \, \d \bX = 0.
\eeq
Applying the divergence theorem to the first term and using the integral constraint, we find
\begin{equation}
    \int_{\partial D_\ee} \bQ_0 \cdot \bn_1 \d S_\bX + \int_{\partial D_\ee} \bX \cdot \sn \bQ_0 \cdot \bn_0 \d S_\bX+ \int_{D_\ee} \sn \cdot \bQ_0 \d \bX  =  -\int_{D_\ee} \rho \d \bX - \int_{D_\ii} \rho \d \bX. \label{eq:wrk1}
\end{equation}
Applying the flux transport theorem to the second term and Reynolds transport theorem to the final term of \eqref{eq:wrk1}, we obtain
\begin{equation}
  \int_{\partial D_e} \bQ_0 \cdot \bn_1  \, \d S_\bX   +\sn \cdot (\varepsilon_{\rm{eff}} \sn \phi_0)   =  - \rho_{\rm{eff}},
\end{equation}
after simplification, where we have defined the effective charge as in \eqref{eq:rhoeff}. Thus, we have the same homogenised problem as obtained with standard approaches \eqref{eq:homprob}, save the presence of an additional boundary term.

 \vspace{0.5em}
\textit{ii.) The correct approach}
\vspace{0.5em}
\begin{figure}
    \centering
    \includegraphics[width = 0.6\textwidth]{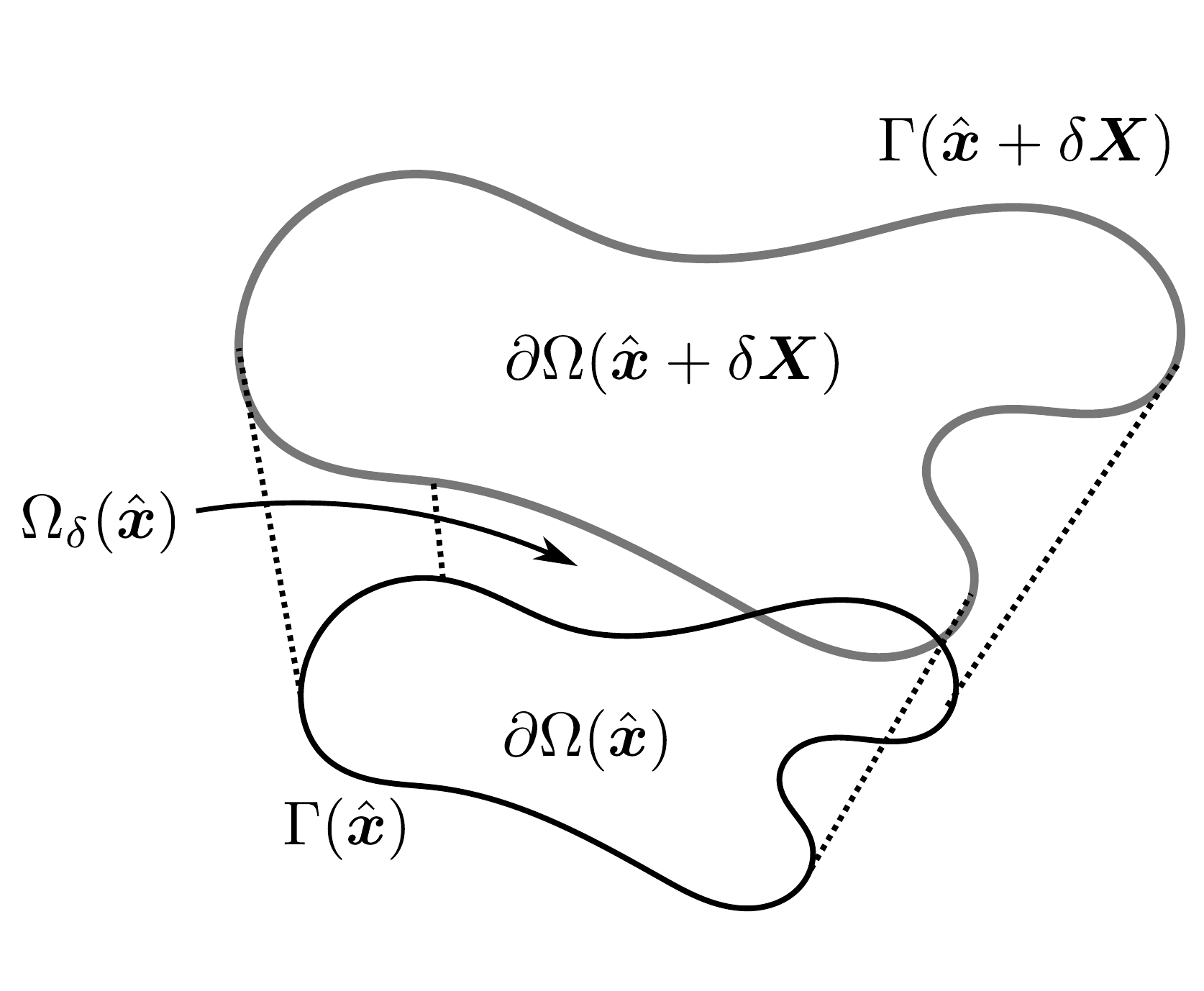}
    \caption{Schematic of the open surface $\partial \Omega$ changing with slow coordinate. }
    \label{fig:flux}
\end{figure}

\noindent In identifying the multiple scales form of integral constraints on a periodic domain, the integrand was expanded about fixed slow position as in \eqref{eq:fastslow}. When the microstructure slowly varies, we must also expand the boundary position about a fixed slow position, \textit{i.e.} we look to expand 
\begin{equation}
    \int_{\partial \Omega} \boldsymbol{Q} \cdot \bn \, \d S =   \delta^2 \int_{\partial \Omega(\hat{\bx} + \delta \bX)} \boldsymbol{Q} (\hat{\bx} + \delta \bX , \bX )\cdot \bn \, \d S_\bX \label{eq:int0} 
\end{equation}
For generality, we assume that the surface $\partial \Omega$ is open as shown in figure \ref{fig:flux}. Expanding the integrand as in \cite{Chapman2015}, we find
\begin{equation}
      \int_{\partial \Omega} \boldsymbol{Q} \cdot \bn \, \d S =   \delta^2 \int_{\partial \Omega(\hat{\bx} + \delta \bX)} \big{(} \boldsymbol{Q} (\hat{\bx} , \bX ) +  \delta \bX \cdot \sn \boldsymbol{Q} (\hat{\bx}  , \bX ))
    \big) \cdot \bn \, \d S_\bX  + \cdots . \label{eq:intxhat}
\end{equation}
To project the boundary onto that at $\hat{\bx}$, we take a similar approach to heuristic derivations of the flux transport theorem, see for example \cite{davis1979}. We apply the divergence theorem to the volume $\Omega_{\delta}$ swept out by the surface as we move from $\hat{\bx}$ to $\hat{\bx} + \delta \bX$ (illustrated schematically in figure \ref{fig:flux}), writing
\begin{equation}
\begin{split}
      \int_{\partial \Omega(\hat{\bx} + \delta \bX)} \boldsymbol{Q} (\hat{\bx} , \bX )  \cdot \bn \,  \d S_\bX  = \int_{\Omega_\delta(\hat{\bx})} \fn \cdot & \boldsymbol{Q}(\hat{\bx} , \bX ) \d \bX
      + \int_{\partial \Omega(\hat{\bx})} \boldsymbol{Q}(\hat{\bx} , \bX ) \cdot \bn \, \d S_\bX  \\&- \int_{\partial \Omega_\delta (\hat{\bx})} \boldsymbol{Q}(\hat{\bx} , \bX )\cdot \bn \, \d S_\bX\label{eq:swept},
\end{split}
\end{equation}
where $\partial \Omega_\delta$ is the volume enclosed by the dashed lines in figure \ref{fig:flux}. In the limit $\delta \to 0$, we can write the volume and surface elements as $\d \bX = \delta \bX \cdot \boldsymbol{V} \cdot \bn \, dS_X$ and $\bn \d S_\bX = -\delta \bX \cdot \boldsymbol{V} \times d \boldsymbol{r}$ respectively, where $\boldsymbol{V} = \sn \boldsymbol{R}^b$ is the `velocity' of the boundary for positions $\boldsymbol{R}^b$ on $\partial \Omega(\hat{\bx})$ and $\d \boldsymbol{r}$ is the line element of $\Gamma(\hat{\bx})$. Substituting the form for the volume element of $\Omega_\delta$ and area element of $\partial \Omega_\delta$ into \eqref{eq:swept}, we have
\begin{equation}
\begin{split}
      \int_{\partial \Omega(\hat{\bx} + \delta \bX)} \boldsymbol{Q} (\hat{\bx} , \bX )  \cdot \bn \,  \d S_\bX  &= \int_{\partial \Omega(\hat{\bx})}  \delta (\fn \cdot  \boldsymbol{Q} ) \bX \cdot \boldsymbol{V} \cdot \bn \,  \d S_\bX  
   \\   +& \int_{\partial \Omega(\hat{\bx})} \boldsymbol{Q} \cdot \bn \, \d S_\bX + \int_{\Gamma (\hat{\bx}) } \delta \boldsymbol{Q}\cdot (\bX \cdot \boldsymbol{V}) \times \d \boldsymbol{r},  \label{eq:xhatexp}
\end{split}
\end{equation}
Combining \eqref{eq:xhatexp} with \eqref{eq:intxhat}, we obtain the multiple scales form for integral constraints on a slowly varying domain
 \begin{equation}
 \begin{split}
     \int_{\partial \Omega}  \bQ \cdot \bn\, \d S  \rightarrow  \delta^2 \int_{\partial \Omega}  \left(\bQ + \delta \bX \cdot \sn  \bQ + \delta (\fn \cdot \bQ) \bX \cdot \boldsymbol{V} \right)\cdot \bn \, \d S_\bX \\
     +  \delta^2 \int_{ \Gamma  } \delta \boldsymbol{Q}\cdot (\bX \cdot \boldsymbol{V}) \times \d \boldsymbol{r} \label{eq:maineqnfull}.
 \end{split}
 \end{equation}
 In \eqref{eq:maineqnfull}, $\bn$ is the normal to the boundary at fixed $\hat{\bx}$; in the example of inclusions with a slowly varying radius this is given by $\bn_0$. Note that, unlike in Section \ref{sec:2dstandard} when approximating \eqref{eq:2decont}, and perhaps counter-intuitively, we do not need to  expand the normal to introduce $\bn_1$, or to apply the operator $\bX \cdot \sn$ to $\bn_0$: the perturbation to the normal is already accounted for by the term involving $\bV$. Thus, an expansion of the normal will only appear in \eqref{eq:maineqnfull} when the function defining the boundary through its level sets is a function of $\delta$.

Thus, in multiple scales form, the integral constraint \eqref{lim:int} is
\beq
\int_{\dd D_\ee}  \left(\bQ + \delta \bX \cdot \sn \bQ + \delta (\fn \cdot \bQ) \bX \cdot \bV\right)\cdot \bn_0\, \d S_\bX = -\delta \int_{D_\ii} \rho \, \d \bX,\label{MSint}
\eeq
where
\[ \bQ = \frac{1}{\delta} \fn \phi + \sn \phi\]
and the surface integral is over the exterior surface of the inclusion.

Using the expansion (\ref{eq:2dphimsexp}) in (\ref{MSint}) we find at leading-order that
\beq
\int_{\dd D_\ee}\fn \phi_0 \cdot \bn_0\, \d S_\bX =0,
\eeq
consistent with $\phi_0 = \phi_0(\bx)$. At first-order we find
\beq
\int_{\dd D_\ee} \bQ_0 \cdot  \bn_0\, \d S_\bX =0,\qquad \qquad \bQ_0 = \fn \phi_1 + \sn \phi_0,
\eeq
which is consistent with \eqref{eq:2dodgauss}.
Finally, equating coefficients of $\delta^2$, and noting that \mbox{$\fn \cdot \bQ_0 = 0$},  we obtain
\beq
\int_{\dd D_\ee}  \left(\bQ_1 +  \bX \cdot \sn \bQ_0 \right)\cdot \bn_0\, \d S_\bX = - \int_{D_\ii} \rho \, \d \bX,\qquad\qquad\bQ_1 =  \fn \phi_2 + \sn \phi_1.
\eeq
Substituting into (\ref{2.49}) gives
\begin{equation}
    \int_{\partial D_\ee} \varepsilon_\ee \bX \cdot \sn \bQ_0 \cdot \bn_0 \,\d  S_\bX + \int_{D_\ee} \sn \cdot \left( \varepsilon_\ee \bQ_0\right) \, \d \bX = - \rho_{{\rm eff}},\label{2.56}
\end{equation}
where 
\begin{equation}
    \rho_{{\rm eff}} = \int_{D} \rho \, \d \bX
\end{equation}
is the effective charge as before. Using the transport theorem to take the slow derivatives outside the integral gives
\begin{eqnarray*}
\int_{\dd D_\ii} (\bX\cdot \sn) \bQ_0\cdot \bn_0\, \d S_\bX & =&
\int_{\dd D_\ii} X_i \frac{\dd Q_{0,j}}{\dd x_i}n_{0,j}\, \d S_\bX\\
& = & \frac{\dd}{\dd x_i} \int_{\dd D_\ii} X_i Q_{0,j}n_{0,j}\, \d S_\bX -  \int_{\dd D_\ii}  Q_{0,i}V_{ij}n_{0,j}\, \d S_\bX\\
& = & \sn \cdot \int_{\dd D_\ii} \bX (\bQ_{0} \cdot \bn_{0})\, \d S_\bX -  \int_{\dd D_\ii}  \bQ_0 \cdot {\bf V} \cdot \bn_{0}\, \d S_\bX,
\end{eqnarray*}
while
\begin{eqnarray*}
  \int_{D_\ee} \sn \cdot \left( \varepsilon_\ee \bQ_0\right) \, \d \bX & = & \sn \cdot\int_{D_\ee}  \varepsilon_\ee \bQ_0 \, \d \bX + \int_{\dd D_\ii }  \bQ_0 \cdot {\bf V} \cdot \bn_{0}\, \d S_\bX
\end{eqnarray*}
(since $\bn_0$ is the outward normal to $D_\ii$). Thus the two surface integrals cancel. Simplifying now as we did to obtain (\ref{eq:permlim}) we find that \eqref{2.56} becomes
\begin{equation}
    \sn \cdotp \left( {\boldsymbol \varepsilon}_{{\rm eff}} \sn \phi_0 \right) = -\rho_{{\rm eff}}, \label{eq:homint}
\end{equation}
where the effective permittivity tensor is
\begin{equation}
   {\varepsilon_{{\rm eff}}}_{ij} = \varepsilon_\ee \left( \delta_{ij} + \int_{\partial D} X_i \frac{\partial \Psi_j}{\partial X_k} n_k \, \d S_\bX  \right) \label{eq:permint}
\end{equation}
in agreement with \eqref{eq:permlim}. 

\section{Paradigm Problem - Another Limit}\label{sec:bigrho}
In section \ref{sec:2d} we illustrated how to treat the multiple scales problem with integral constraints considered in \cite{Chapman2015} when the domain slowly varies. In this example, the divergence of the flux in the proposed multiple scales form \eqref{MSint} vanishes. Here we construct a problem where this term is non-zero by considering the limit of large charge density, rescaling \eqref{eq:2dgauss}-\eqref{lim:int} as follows. We consider
\begin{equation}
    \nabla \cdot ( \varepsilon \nabla \phi) = -\frac{\rho}{\delta} \label{eq:2dgaussb},
\end{equation}
where $\rho = O(1)$. The boundary conditions \eqref{eq:2decont} and \eqref{eq:2dcont} remain unchanged. 

In the limit of perfectly dielectric inclusions, where $\varepsilon_\ii \to \infty$, we have
\begin{align}
  \nabla \cdot ( \varepsilon_\ee \nabla \phi) & = -\frac{\rho}{\delta}\qquad \mbox{ in }\Omega_\ee,\label{lim:eqeb}\\
  \nabla \phi & = 0\qquad \mbox{ in }\Omega_\ii,\label{lim:eqib}
\end{align}
with continuity \eqref{lim:2dcont} at the inclusion boundary. 
The rescaled integral constraint becomes
\beq
\int_{\dd \Omega_\ii} \left.\varepsilon_\ee\, \bn\cdot \del \phi \right|_\ee\, \d S = - \frac{1}{\delta} \int_{\Omega_\ii} \rho \, \d \bx. \label{lim:intb}
\eeq
We perform a similar analysis to section \ref{sec:2d}: first we take the limit $\varepsilon_\ii \to \infty$ in the standard multiple scales problem before comparing with the problem formulated with an integral condition.

\subsection{Standard Multiple Scales}\label{sec:2dstandardb}
We substitute the multiple scales expansion \eqref{eq:2dphimsexp} into \eqref{eq:2dgaussb} and compare coefficients at each order of $\delta$. At leading-order, we find $\phi_0$ independent of the fast scale. At next order, we have
\begin{align}
    \fn \cdot \left( \varepsilon \fn \phi_1 \right) &= -\rho \quad \text{in} \quad D, \\
     \left[ \bn_0 \cdot \left( \varepsilon(  \fn \phi_1 +  \sn  \phi_0)  \right)  \right]^\ee_\ii &= 0,\label{eq:odb} \\
 \left[ \phi_1 \right]^\ee_\ii &= 0,
\end{align}
with $\phi_1$ $\bf 1$-periodic in $\bf X.$ Writing 
\begin{equation}
    \phi_1 = {\bf \Psi}\cdot \sn \phi_0 + \xi +  \overline{\phi}_1,\label{eq:2dphi1b}
\end{equation}
with $\overline{\phi}_1 = \overline{\phi}_1({\bf x})$, we obtain two microscale problems. We find ${\bf \Psi}$ satisfies \eqref{eq:potcell2d}-\eqref{eq:potcellcont}.
The microscale problem for $\xi$ is
\begin{align}
     \fn \cdot \left( \varepsilon \fn {\xi} \right) &= -\rho \quad \text{in} \quad D,  \label{eq:potcell2dxi}\\
    \left[ \bn_0  \cdot\varepsilon \fn {\xi} \right]^\ee_\ii &= 0, \label{eq:neumannxi}\\
    \left[ {\xi} \right]^\ee_\ii &= 0,\label{eq:potcellcontxi}
\end{align}
with 
\begin{equation}
    \int_{D} \xi \, d\bX = 0. 
\end{equation}
Equating coefficients of $\delta^2$, we find
\begin{align}
  \fn \cdot \left( \varepsilon ( \fn \phi_2  + \sn \phi_1   )\right) + \sn \cdot \left( \varepsilon ( \fn \phi_1  + \sn \phi_0   )\right) &= 0  \quad \text{in} \quad D,  \label{eq:2dod2gaussb} \\
     \left[ \bn_0 \cdot \left( \varepsilon(  \fn \phi_2 +  \sn  \phi_1)  \right) \right]^\ee_\ii  +   \left[ \bn_1 \cdot \left( \varepsilon(  \fn \phi_1 +  \sn  \phi_0)  \right)  \right]^\ee_\ii &= 0, \label{eq:2dod2jumpb}\\
 \left[ \phi_2 \right]^\ee_\ii &= 0.\label{eq:2dod2contb}
\end{align}
We integrate (\ref{eq:2dod2gaussb}) over the unit cell $D$, applying the divergence theorem to terms involving the fast divergence, the Reynolds transport theorem and use (\ref{eq:2dod2jumpb}). Following similar analysis to section \ref{sec:2dstandard}, we obtain the homogenised problem
\begin{equation}
  \sn \cdot \left( {\boldsymbol \varepsilon}_{{\rm eff}}  \sn \phi_0 \right)
  = -\rho_{\rm{eff}}, \label{eq:homb}
\end{equation}
with effective permittivity ${\boldsymbol \varepsilon}_{{\rm eff}} $ given by
\eqref{2.31} and effective charge density,
\begin{equation}
    \rho_{\rm{eff}} = \sn \cdotp \int_{D} \varepsilon \fn \xi d\bf X.\label{eq:rhoeffb}
\end{equation}
\subsubsection{Taking the limit $\varepsilon_{\ii} \to \infty$}
In the limit of perfectly dielectric inclusions, the effective permittivity takes the form \eqref{eq:permlim}. In the limit $\varepsilon_{\ii} \to \infty$, the cell problem for $\xi $ becomes
\begin{align}
    \fn \cdot (\varepsilon_\ee \fn \xi) &= - \rho \quad \text{in} \quad D_\ee,\label{eq:xiext}\\
    \fn^2 \xi &= 0 \quad \text{in} \quad D_\ii,\label{eq:xiint}\\
     \bn_0 \cdotp \fn \xi &= 0 \quad \text{on} \quad \partial D_\ii, \\
     \left[ \xi  \right]^\ee_\ii &= 0 .
\end{align}
We switch to index notation to establish the form of the effective charge $\rho_{\rm{eff}}$ in the limit $\varepsilon_{\ii} \to \infty$, 
\begin{equation}
\begin{split}
    \rho_{\rm{eff}} &= \frac{\dd }{\dd x_i} \int_{D} \varepsilon \frac{\dd \xi}{\dd X_i} \d {\bf X} \\
    &= \frac{\dd }{\dd x_i} \int_{D} \varepsilon \frac{\dd }{\dd X_j}\left( X_i \frac{\dd \xi}{\dd X_j} \right) \d {\bf X} - \frac{\dd}{\dd x_i} \int_{D} \varepsilon X_i \frac{\dd^2 \xi }{\dd X_j \dd X_j} \d {\bf X}\\
   &=  \frac{\dd }{\dd x_i} \int_{\partial D} \varepsilon
   _\ee  X_i \frac{\dd \xi}{\dd X_j}  n_j \d S - \frac{\dd}{\dd x_i} \int_{D} \varepsilon X_i \frac{\dd^2 \xi }{\dd X_j \dd X_j} \d {\bf X}\\
     &= \frac{\dd }{\dd x_i} \int_{\partial D} \varepsilon_\ee X_i \frac{\dd \xi}{\dd X_j} n_j \d S + \frac{\dd}{\dd x_i} \int_{D_e}  X_i \rho \d {\bf X},
\end{split}
\end{equation}
where we have used \eqref{eq:neumannxi} in going from the second to third line and \eqref{eq:xiext}-\eqref{eq:xiext} in the third to fourth line.

\subsection{Multiple Scales with Integral Constraints}
In this section we treat \eqref{lim:eqeb}-\eqref{lim:intb} directly, writing \eqref{lim:intb} in multiple scales form using \eqref{eq:maineqnfull}. Substituting \eqref{eq:2dphimsexp} into \eqref{lim:eqeb}-\eqref{lim:intb} written in multiple scales form, we find that $\phi_0 = \phi_0({\bf x})$ at leading-order.  At first-order we find
\begin{align}
    \fn\cdot ( \varepsilon_\ee \fn \phi_1 ) & = -\rho \quad \text{in} \quad D_\ee, \label{eq:2dodgaussb} \\
    \fn \phi_1 + \sn \phi_0 &= {\bf 0} \quad \text{in} \quad D_\ii,\label{eq:odincb}\\
    \left[ \phi_1  \right]^\ee_\ii &= 0, \\
    \int_{\partial D_\ee} \varepsilon_\ee (\fn \phi_1 + \sn \phi_0 ) \cdotp \bn_0 \d S_\bX  & = - \int_{D_\ii} \rho \d {\bf X},
\end{align}
with $\phi_1$ ${\bf 1}$-periodic in $\bX$. As in Section \ref{sec:2dstandardb}, we write $\phi_1 = {\bf \Psi}\cdot \sn \phi_0 + \xi +  \overline{\phi}_1$ where
$\overline{\phi}_1 = \overline{\phi}_1({\bf x})$ and ${\bf \Psi}$ is the solution to \eqref{eq:cellint}-\eqref{eq:cellintcont}. 
The second cell function $\xi$ satisfies
\begin{align}
        \fn \cdot ( \varepsilon_\ee \fn \xi ) &= -\rho \quad \text{in} \quad D_\ee, \label{eq:xiinf} \\
    \fn \xi &= {\bf 0} \quad \text{in} \quad D_\ii, \\
    \left[ \xi \right]^\ee_\ii &= 0, \label{eq:xiinfcont}
\end{align}
with 
\begin{equation}
    \int_{D} \xi \d{\bf X  } = 0.
\end{equation}
Equating coefficients at next order, we have
\begin{align}
    \fn \cdot \left( \varepsilon_\ee (\fn \phi_2 + \sn \phi_1 )\right) +   \sn \cdot \left( \varepsilon_\ee (\fn \phi_1 + \sn \phi_0 )\right) &= 0 \quad \text{in} \quad D_\ee, \label{eq:od2gaussintb}\\
    \fn \phi_2 + \sn \phi_1 &= {\bf 0} \quad \text{in} \quad D_\ii, \\
    \left[ \phi_2 \right]^\ee_\ii &= 0,\\
    \int_{\partial D_\ee} \varepsilon_\ee (\fn \phi_2 + \sn \phi_1 ) \cdotp \bn_0 \d S_\bX +   \int_{\partial D_\ee} \varepsilon_\ee  \bX\cdotp \sn  (\fn \phi_1 +& \sn \phi_0 )  \cdotp \bn_0 \d S_\bX \nonumber\\+   \int_{\partial D_\ee} \varepsilon_\ee  \fn \cdotp(\fn \phi_1 + \sn \phi_0 ) {\bf X \cdotp \bf  V }\cdotp \bn_0 \d S_\bX   &= 0\label{eq:od2intb}. 
\end{align}
Integrating (\ref{eq:od2gaussint}) over the exterior region, applying the divergence theorem to the first term and substituting \eqref{eq:od2intb}, we have
\begin{equation}
\begin{split}
 \int_{\partial D_\ee} \varepsilon_\ee  \bX\cdotp \sn  (\fn \phi_1 &+ \sn \phi_0 ) \cdotp \bn_0 \d S_\bX +   \int_{\partial D_\ee} \varepsilon_\ee  \fn \cdotp(\fn \phi_1 + \sn \phi_0 ) {\bf X \cdotp \bf  V }\cdotp \bn_0 \d S_\bX  \nonumber \\ &+ \int_{D_\ee} \sn \cdot \left( \varepsilon_\ee (\fn \phi_1 + \sn \phi_0 )\right)\, \d \bX = 0.
 \end{split}
\end{equation}
Applying the transport theorem to the first and final integrals gives
\begin{equation}
\begin{split}
 \sn \cdotp \int_{\partial D_\ee} \varepsilon_\ee  \bX  (\fn \phi_1 &+ \sn \phi_0 ) \cdotp \bn_0 \d S_\bX -   \int_{\partial D_\ee} \varepsilon_\ee \fn \cdotp ( \bX  (\fn \phi_1 + \sn \phi_0 ) ) {\cdotp \bf V} \cdotp \bn_0 \d S_\bX \\ 
 +   \int_{\partial D_\ee} \varepsilon_\ee  \fn \cdotp(\fn \phi_1 &+ \sn \phi_0 ) {\bf X \cdotp \bf  V }\cdotp \bn_0 \d S_\bX  + 
 \sn \cdot \int_{D_\ee}  \left( \varepsilon_\ee (\fn \phi_1 + \sn \phi_0 )\right)\, \d \bX \\ &+    \int_{\partial D_\ee} \varepsilon_\ee (\fn \phi_1 + \sn \phi_0 )  {\cdotp \bf V} \cdotp \bn_0 \d S_\bX   = 0.\nonumber
\end{split}
\end{equation}
Expanding the divergence in the second integral, we find some of the boundary terms cancel, leaving
\begin{equation}
    \sn \cdotp \int_{\partial D_\ee} \varepsilon_\ee  \bX  (\fn \phi_1 + \sn \phi_0 ) \cdotp \bn_0 \d S_\bX +  \sn \cdot \int_{D_\ee}  \left( \varepsilon_\ee (\fn \phi_1 + \sn \phi_0 )\right)\, \d \bX = 0.
\end{equation}
We substitute \eqref{eq:2dphi1b}, using the divergence theorem to take the first integral into the exterior region and use \eqref{eq:xiinf}-\eqref{eq:xiinfcont} to obtain
\begin{equation}
    \sn \cdot (\varepsilon_{\rm{eff}} \sn \phi_0 ) = - \rho_{\rm{eff}}
\end{equation}
where 
\begin{equation}
    \varepsilon_{\rm{eff} ij} = \varepsilon_{\ee} \left( \delta_{ij} + \int_{\partial D} X_i \frac{\dd \Psi_j}{\dd X_k } \bn_{0k} \d S_\bX \right).
\end{equation}
The effective charge is given by
\begin{equation}
\rho_{\rm{eff}} = \sn\cdotp \left( \int_{D_e}\rho {\bf X} \d {\bf X} + \int_{\partial D} \varepsilon_\ee {\bf X} \fn \xi \cdot \bn_0  \d S_\bX \right).
\end{equation}
Thus, we have recovered \eqref{eq:homb} in the limit of perfectly dielectric inclusions, confirming the need for the divergence term present in \eqref{MSint}.

\section{Discussion}
We have outlined how to combine the extension to the standard theory of multiple scales which deals with a slowly varying microstructure with that which deals with  integral constraints.
Our main result is equation (\ref{eq:maineqnfull}), which shows how to write an integral constraint in multiple scales form when the (fast) domain of the integral is a function of the slow scale. Essentially the rest of the manuscript is a justification of this equation, showing that it leads to the correct homogenised model for an example in which that model can be identified using a more standard approach.
Some problems involving integral constraints, especially those in which different physics holds in the inclusions, do not arise as a limit of a more standard problem, and such an approach is not available.
These problems can be handled using equation (\ref{eq:maineqnfull}).

\section{Acknowledgements}
A. K. thanks the BBSRC for support under grant BB/M011224/1.

\label{lastpage}

\end{document}